\documentclass[12pt]{amsart}
\usepackage{amssymb,amsthm}

\def\silentfootnote#1{{\let\thefootnote\relax\footnotetext{#1}}}

\numberwithin{equation}{section}
\theoremstyle{plain}
\newtheorem{thm}[equation]{Theorem}
\newtheorem{lem}[equation]{Lemma}
\newtheorem{prop}[equation]{Proposition}

\newtheorem{ex}[equation]{Example}
\newtheorem{rem}[equation]{Remark}

\newcommand{\Z}{{\mathbb Z}}
\newcommand{\ra}{\rightarrow}
\newcommand{\C}{{\mathbb C}}

\DeclareMathOperator{\Hom}{Hom} 
\DeclareMathOperator{\Ext}{Ext} \DeclareMathOperator{\Ker}{Ker}
\DeclareMathOperator{\Ima}{Im} 
 \DeclareMathOperator{\coh}{H}
\DeclareMathOperator{\HH}{HH} \DeclareMathOperator{\id}{id}

\DeclareFontFamily{U}{rsf}{} \DeclareFontShape{U}{rsf}{m}{n}{
  <5> <6> rsfs5 <7> <8> <9> rsfs7 <10-> rsfs10}{}
\DeclareMathAlphabet{\mathscr}{U}{rsf}{m}{n}

\DeclareMathAlphabet{\mathgth}{U}{euf}{m}{n}

\DeclareFontFamily{U}{cyr}{} \DeclareFontShape{U}{cyr}{m}{n}{
  <5> wncyr5 <6> wncyr6 <7> wncyr7 <8> wncyr8 <9> wncyr9 <10-> wncyr10}{}
\DeclareMathAlphabet{\mathcyr}{U}{cyr}{m}{n}

\input cyracc.def

\newcommand{\cO}{{\mathscr O}}

\newcommand{\D}{{\mathbf D}_{\mathrm{coh}}^b}

\DeclareMathOperator{\Spec}{Spec}

\DeclareMathOperator{\Proj}{Proj}

\DeclareMathOperator{\Physics}{Physics}
\DeclareMathOperator{\Def}{Def}

\newcommand{\iso}{\cong}

\renewcommand{\phi}{\varphi}

\title[Discrete Torsion]
{Algebraic Deformations Arising from Orbifolds with Discrete
Torsion}
\author{Andrei C\u ald\u araru}
\address{Department of Mathematics\\
         University of Pennsylvania\\
         Philadelphia, Pennsylvania 19104}
\email{andreic@math.upenn.edu}

\author{Anthony Giaquinto}
\address{Department of Mathematics and Statistics\\
      Loyola University Chicago\\
     Chicago, Illinois 60626}
\email{tonyg@math.luc.edu}

\author{Sarah Witherspoon}
\address{Department of Mathematics and Computer Science\\
         Amherst College\\
         Amherst, Massachusetts 01002}
\email{sjw@cs.amherst.edu}

\date{June 6, 2003}

\begin{document}

\begin{abstract}
We develop methods for computing Hochschild cohomology groups and
deformations of crossed product rings.  We use these methods to find
deformations of a ring associated to a particular orbifold with
discrete torsion, and give a presentation of the center of the
resulting deformed ring.  This connects with earlier calculations by
Vafa and Witten of chiral numbers and deformations of a similar orbifold.
\end{abstract}

\maketitle

\silentfootnote{The third author was supported by National
Security Agency Grant \#MDS904-01-1-0067.}



\section{Introduction}

Our motivation for this paper was a desire to provide a mathematical
basis for physics statements in the paper \cite{VafWit} of Vafa and
Witten.  Specifically, we expect that the chiral numbers of an
orbifold with discrete torsion come from the Hochschild cohomology
groups of associated crossed product rings, and that geometric
deformations of the orbifold correspond with algebraic deformations of
these rings.  In Section 2 we make these ideas more precise, and
propose a mathematical definition of some of the chiral numbers.  This
gives the proper context for results in the remainder of the paper.
We also describe in detail there the central example of \cite{VafWit},
which consists of the quotient of the product of three elliptic curves
by an action of the Klein four-group.  A local version of this is our
motivating example throughout the paper.

Specifically, let $R=\C [x,y,z]$, $X=\C^3=\Spec R$, and let
$G=\Z/2\times \Z/2$ act on $X$ by pairwise negation on two
coordinates, leaving the third one fixed.  There are three curves in
$X$ on which the action of $G$ is not free (the coordinate axes), all
meeting in the origin, which is fixed by $G$.  Details are given in
Example \ref{example}.  We also include in Section 3 the definitions
(more generally) of the crossed product ring $R\#_{\alpha}G$ (where
$\alpha$ is a two-cocycle), and of the Hochschild cohomology groups
and deformations of a ring.

Results of \c Stefan \cite{stefan95} imply that (over $\C$)
the Hochschild cohomology group $\HH^*(R\#_{\alpha}G)$
is the subspace of $\HH^*(R,R\#_{\alpha}G)$ invariant under an action of $G$.
Letting $R=\C[x,y,z]$, we will use the Koszul resolution of $R$ over $R\otimes R$ to
find $\HH^*(R,R\#_{\alpha}G)$ explicitly, and we will want an action
of $G$ on the Koszul complex.
This action is given in Lemma \ref{action} in case $G$ is a finite
abelian group.
We apply these ideas to find $\HH^*(R\#_{\alpha}G)$ when
$G=\Z/2\times \Z/2$ in Example \ref{hh*}.
For this example, unlike the compact example of Vafa and Witten,
the Hochschild cohomology groups that we compute are infinite dimensional
as vector spaces over $\C$.

In general, elements of $\HH^2(R\#_{\alpha}G)$ coincide with infinitesimal
(or first order) deformations of the multiplication of $R\#_{\alpha}G$.
In order to find {\em explicit} formulas for the deformations in case
$R=\C[x,y,z]$, we must relate the Koszul complex
for $R$ to the Hochschild complex for $R\#_{\alpha}G$.
In Section 5 we develop a general method for doing so, and apply it
to find the infinitesimal deformations of our example
$A=\C[x,y,z]\#_{\alpha}(\Z/2\times \Z/2)$.
In Section 6 we continue and show that for our example, in the
presence of discrete torsion ($\alpha$ nontrivial),
the infinitesimal deformations of $R\#_{\alpha}G$
not coming from $R$ lift to formal deformations of $R\#_{\alpha}G$.
This follows from the application of a universal deformation
formula based on a particular bialgebra, using results of 
\cite{giaquinto-zhang98}.
Our example is of independent interest in this context:
it is the first example that we know of a formal deformation arising
from a {\em noncocommutative} bialgebra.

In the last section of the paper we study the behavior of the center
$Z(R\#_{\alpha}G)$, which is simply the $G$-invariant subring $R^G$,
under deformation.  In our example,
the three curves of singularities of $\Spec R^G$ are smoothed out
after formal deformation, and an isolated singularity is left at the
origin.
We give explicit equations for this singularity in terms of the
polynomials arising in the first order deformation.

We have now described briefly some of the calculations in this paper.
They indicate that the behavior of our local picture
closely mimics the picture described by Vafa and Witten in
their example.
This supports our belief that the mathematical study of an
orbifold with discrete torsion should begin with the study of the
Hochschild cohomology groups and deformations of associated crossed
product rings, followed by an understanding of the behavior
of the centers of these rings under deformation.

\section{Physical and mathematical context}

We begin with some motivating ideas from physics and describe
in particular an example of Vafa and Witten.

\subsection{Chiral numbers for Calabi-Yau manifolds and orbifolds}
Given an action of a finite group $G$ on a space $X$ satisfying
certain properties, one can construct several physics theories
which represent the $G$-equivariant
physics of $X$.  There is more than one such theory because their
construction takes as additional input a two-cocycle
$\alpha:G\times G\rightarrow \C^{\times}$, called the
{\em discrete torsion} of the theory.  We denote (in a loose
sense) the $G$-equivariant physics theory on $X$, with discrete
torsion $\alpha$, by $\Physics^{G,\alpha}(X)$.  A similar
construction exists in the non-equivariant setting (without any
contribution from discrete torsion), assigning to a Calabi-Yau
threefold $X$ a physics theory $\Physics(X)$.

Every such theory has associated to it {\em chiral numbers}
$h^{ij}$ ($i,j=0,\ldots, 3$) which are defined as certain physical
quantities.  We are interested in giving a mathematical definition
of these numbers in situations of interest to physicists.

For a non-equivariant theory arising from an ordinary
Calabi-Yau threefold $X$ it is well known that
\[ h^{ij}(\Physics(X)) = h^{ij}(X), \]
the ordinary Hodge numbers of $X$.  An important observation is
that $h^{12}(X)$ is the dimension of the space of {\em
complex deformations} of $X$, which is also the
dimension of the space of deformations of the multiplication of
the structure sheaf $\cO_X$ of $X$.

Chiral numbers are also well understood for a theory
arising from a group action $G$ on a space $X$, but with trivial
discrete torsion.  Physics predicts that the chiral
numbers of such a theory should be the Hodge numbers of a crepant
resolution of the singularities of $X/G$.\footnote{A resolution
$Y\ra X/G$ of the singularities of $X/G$ is {\em
crepant} if the relative canonical bundle is trivial.  In situations of interest in physics, this is equivalent to $Y$ being Calabi-Yau.}
This follows from the general principle
of Kontsevich that (at least some of) the data contained in the
physics model can be extracted from a certain derived category.
For a theory built from an ordinary Calabi-Yau
threefold $X$ this category is $\D(X)$, the derived category of
coherent sheaves on $X$. In the equivariant setting this should be
replaced by the derived category $\D([X/G])$ of $G$-{\em equivariant}
sheaves on $X$.  Results of Bridgeland, King and Reid~\cite{BKR},
combined with further results of Bridgeland~\cite{BriFlops} show
that in the cases of interest in physics there is an equivalence
\[ \D([X/G]) \iso \D(Y) \]
for any crepant resolution $Y$ of the singularities of $X/G$.
Since one should be able to extract $h^{ij}(\Physics^G(X))$ from
$\D([X/G])$, the above isomorphism would yield
\begin{align*}
h^{ij}(\Physics^G(X)) & = h^{ij}(\D([X/G])) = h^{ij}(\D(Y)) =
h^{ij}(\Physics(Y)) \\
& = h^{ij}(Y).
\end{align*}

This argument is unsatisfactory for two reasons: first, there is no good
notion of Hodge numbers for an arbitrary derived category, so the
above equalities should be thought of only as general principles,
not as mathematical statements.  Second, there is no generalization of
the above argument regarding chiral numbers
in the presence of discrete torsion.

It should be mentioned here that there is another, topological
approach to chiral numbers for orbifolds, called {\em orbifold
cohomology}. See~\cite{Chen-Ruan} and~\cite{Lupercio-Uribe} for details.

\subsection{Hochschild cohomology}
Noncommutative geometry provides a different ap\-proach.
Some of the following ideas have been around for a
while, see for example Connes' book~\cite{Con}.

If $R$ is a commutative ring with an action of a finite group $G$,
we can construct the crossed product ring $R\#G$ (see
Section~\ref{sec:prel} for a definition).  It is an
associative, noncommutative ring with the property that
$R\#G$-modules correspond with $G$-equivariant
$R$-modules. When the finite group $G$ acts on a
scheme $X$, this construction
yields a coherent sheaf $\cO_X\#G$ of noncommutative
algebras on $X$, with the property that coherent sheaves
of $\cO_X\#G$-modules are identified with
$G$-equivariant coherent sheaves of $\cO_X$-modules on $X$.
More generally if $\alpha\in\coh^2(G,\C^{\times})$, we may construct
a crossed product ring $R\#_{\alpha}G$ twisted by $\alpha$.
Again this construction globalizes to schemes.

The proper choice of cohomology theory for noncommutative rings is
Hoch\-schild cohomology, which is defined for sheaves of algebras in
\cite{gerstenhaber-schack88}.  Hochschild {\em ho}mology is defined
for schemes in \cite{Wei}, and even more generally, Hochschild
homology is defined for exact categories~\cite{Kel,mccarthy94}.
Standard formalism gives as well Hochschild {\em co}homology for exact
categories \cite{mccarthy}.  Thus one can consider the Hochschild
cohomology of the exact category of sheaves of $\cO_X\#_\alpha
G$-modules on a scheme $X$.  Hochschild homology has been shown to be
invariant with respect to derived equivalences coming from
localization of pairs~\cite{Kel}, and the same techniques can be
applied for cohomology as well.  In fact Hochschild cohomology for a
scheme $X$ can be shown to be invariant under Fourier-Mukai
transforms~\cite{Cal} (following ideas in~\cite{Kel1}, where the
affine case is studied).  In particular, the Bridgeland-King-Reid
equivalence
\[ \D([X/G]) \iso \D(Y), \]
to which we referred earlier, yields an isomorphism on Hochschild cohomology
\[ \HH^i(\cO_X\#G) \iso \HH^i(\cO_Y) \]
for every $i$.  For a smooth quasi-projective scheme $Y$ we
have~\cite{gerstenhaber-schack88,Kon}
\[ \HH^i(\cO_Y) \cong \bigoplus_{p+q=i}\coh^p(Y, \bigwedge^q T_Y), \]
where $T_Y$ is the tangent bundle of $Y$.
Thus for a smooth Calabi-Yau threefold $Y$, for which we have
\[ \bigwedge^q T_Y \iso \bigwedge^{3-q} \Omega_Y, \]
the dimensions of $\HH^i(\cO_Y)$ are given by
\[ 1,~0,~h^{12}(Y),~2\cdot h^{11}(Y) + 2,~h^{12}(Y),~0,~1. \]
In the case of a Calabi-Yau threefold then, the physical statement that
\[ h^{ij}(\Physics^G(X)) = h^{ij}(\Physics(Y)) \]
can be viewed as a consequence of the mathematical statement that
\[ \HH^i(\cO_X \# G) \iso \HH^i(\cO_Y) \]
for all $i$.  (Note that knowing the dimensions of the Hochschild
cohomology groups for a Calabi-Yau threefold allows us to recover
its Hodge numbers $h^{11}$ and $h^{12}$.)  In higher dimensions
we do not expect to recover all the chiral numbers but rather only
the sums
\[ \dim \HH^i = \sum_{n+p-q = i} h^{pq}, \]
where $n$ is the dimension of the underlying space.

In the case of a Calabi-Yau threefold $X$, we could define the chiral
numbers $h^{11}$ and $h^{12}$ as
\[ h^{11} = \frac{1}{2}\dim \HH^3(\cO_X) -1 \ \mbox{ and } \ h^{12} =
\dim \HH^2(\cO_X). \] 
The same definition makes sense for an orbifold $[X/G]$ (possibly with
discrete torsion $\alpha$), in which case we replace $\cO_X$ by
$\cO_X\#_\alpha G$.  This allows us to talk about chiral numbers in
the generalized setting of orbifolds.

A further argument that this definition is the right one is the
following observation: for an arbitrary algebra $A$, $\HH^2(A)$ is the
space of infinitesimal (first order) deformations of $A$ (see Section
3).  In the case $Y$ is a Calabi-Yau threefold, we know that these
deformations are measured by $h^{12}(X)$.  Thus defining $h^{12} =
\dim \HH^2$ for Calabi-Yau-like spaces agrees with the interpretation
of these numbers as dimensions of first order deformation spaces.

\subsection{The example of Vafa and Witten~\cite{VafWit}}
Let $X=E_1\times E_2\times E_3$ be the product of three elliptic
curves, let $G=\Z/2\times \Z/2$, and consider the action of $G$ on
$X$ in which every non-identity element of $G$ acts by negation on
two of the $E_i$'s and leaves the third one fixed.  Note that the
nonzero holomorphic 3-form on $X$ is invariant with respect to
the action of $G$, so that $X/G$ is in a natural sense a
(singular) Calabi-Yau 3-space.  There are 48 curves in $X$
where the action is not free, and these curves intersect in the 64 points
that are fixed by all of $G$.

The discrete torsion group (or Schur multiplier)
$\coh^2(G, \C^{\times})$ is equal to $\Z/2$.
Thus there are precisely two physical theories that we can
build on the orbifold $[X/G]$ corresponding to the two choices of
discrete torsion $\alpha\in \coh^2(G,\C^{\times})$: no discrete
torsion ($\alpha=1$) and nontrivial discrete torsion ($\alpha\neq 1$).

Vafa and Witten computed the chiral numbers of these physical
theories \cite{VafWit},
\begin{align*}
h^{11} = 51,~h^{12} = 3&\mbox{ when }\alpha=1, \\
h^{11} = 3,~h^{12} = 51&\mbox{ when }\alpha\neq 1.
\end{align*}
The case $\alpha=1$ is again perfectly understood: essentially
$h^{12} = 3$ means that every deformation of $\cO_X\#G$ must arise
from a deformation of $X$ (note that $h^{12}(X) = 3$ corresponding
to the possible changes of complex structure on each elliptic
curve).  The fact that $h^{11} = 51$ corresponds to the fact that
to obtain a crepant resolution of $X/G$ one needs to blow-up the
48 curves of singularities in $X/G$, thus obtaining 48 new
K\"ahler deformations of the resolution $Y$.  (Again, $h^{11}(X) =
3$, corresponding to a choice of volume for each elliptic curve.)

The surprise lies in the value $h^{12} = 51$ in the presence of
discrete torsion.  In general, given a singular theory whose
singularities are resolved by blow-ups (as in the resolution $Y\ra
X/G$), physics predicts the existence of another theory which removes
the singularities by deforming them.  Thus one {\em does} expect to
have more deformations.  The problem is that the total number of
deformations of $X/G$ is 115, which exceeds $h^{12}$, and so some
deformations of $X/G$ are not allowed in the physical model.  Vafa and
Witten guessed that this means the allowed deformations of $X/G$ only
partially smooth $X/G$, that is they are required (for some mysterious
reason) still to have 64 ordinary double points.  An explanation for
the appearance of these 64 ordinary double points has since been given
by the first author~\cite{CalVW}, following ideas of Aspinwall,
Morrison and Gross~\cite{AspMorGro}.

{F}rom the perspective of noncommutative geometry, this
situation can be described as follows: the data of
$\Physics^{G,\alpha}(X)$ (both with and without discrete torsion)
is encoded in the derived category
$\D(\cO_X\#_\alpha G)$.
Deformations of the ring $\cO_X\#_\alpha G$ are measured by the
chiral number $h^{12}$, and each deformation of this ring yields a
deformation of its center, which coincides with the ring of
invariants
$(\cO_X)^G = \cO_{X/G}$.
In each case we get an allowable set of deformations of
$\cO_{X/G}$. In the absence of discrete torsion, there is only a
3-dimensional space of allowable deformations (that all come from
deformations of $X$), while in the presence of discrete torsion,
there is a 51-dimensional space of deformations, all of which
deform the center $\cO_{X/G}$ to the structure sheaf of a space
which has at least 64 ordinary double points.

\subsection{Our example}
For the purposes of our calculations we simplify Vafa and Witten's
example by replacing the product of three elliptic curves by affine
3-space.  One important difference is that in our case, deformations
of arbitrary degree are allowed.  (The situation is analogous to the
difference between affine and projective geometry: for example if $X$
is the affine scheme $\Spec k[x]$, then $\coh^0(X, \cO_X)$ consists of
all polynomials in $k[x]$ of arbitrary degree, while if $X'$ is the
projective scheme $\Proj k[x,y]$, the global sections $\coh^0(X',
\cO_{X'}(n))$ of a line bundle $\cO_{X'}(n)$ will correspond only to
homogeneous polynomials of fixed degree $n$).  In our example then,
the Hochschild cohomology groups $\HH^*(R\#_{\alpha}G)$ that we will
compute will be infinite dimensional.  When $\alpha =1$, there is
another difference between our local picture and the one described by
Vafa and Witten: $\HH^2(R\#_{\alpha}G)$ picks up many first order
deformations of $A$ that do {\em not} arise from deformations of $R$
(see Example \ref{hh*}).

The deformations of $Z(R\#_{\alpha}G)=R^G$, corresponding to those of
$R\#_{\alpha}G$ that we will calculate, exhibit different behaviors
for different degrees involved.  Similar to the example of Vafa and
Witten, in all cases the three curves of singularities of $\Spec R^G$
are smoothed out after a formal deformation and an isolated
singularity is left at the origin.  However, a surprising consequence
of our calculations is that in order to see this, calculations to
first order do not suffice: first order deformations of
$R\#_{\alpha}G$ do {\em not} give rise to nontrivial {\em first} order
deformations of $R^G$ (see Section 7).

This could be described geometrically as follows: if there were a
geometric germ of a moduli space $(\Def_{R\#_{\alpha}G}, 0)$ of formal
deformations of $R\#_{\alpha}G$ such that maps from $\Spec k[[t]]$ to
$\Def_{R\#_{\alpha}G}$ that map the closed point of $\Spec k[[t]]$ to
$0$ correspond to formal deformations of $R\#_{\alpha}G$, then the
tangent space to $\Def_{R\#_{\alpha}G}$ at $0$ would be naturally
isomorphic to $\HH^2(R\#_{\alpha}G)$, the space of first-order
deformations of $R\#_{\alpha}G$.  The same picture can be imagined for
$R^G$ with a germ of a moduli space $(\Def_{R^G}, 0)$.  The operation
of taking a ring to its center would give a (partially defined) map
$\Def_{R\#_{\alpha}G} \ra \Def_{R^G}$, which in turn would induce a
partially defined map on tangent spaces at the origin
$\HH^2(R\#_{\alpha}G) \ra \HH^2(R^G)$.  (The reason this map is only
partially defined is that there are flat deformations of
$R\#_{\alpha}G$ over $\Spec k[[t]]$ whose center is {\em not} flat
over $\Spec k[[t]]$.)  Our statement in the previous paragraph amounts
to the statement that this map on tangent spaces is zero whenever it
is defined.  In other words, the map $\Def_{R\#_{\alpha}G} \ra
\Def_{R^G}$ is completely ramified at $0$.  However, when we move from
the map on tangent spaces to the map on deformation spaces the map
stops being constant, and this yields the change in the type of
singularity.


\section{Definitions}
\label{sec:prel}

In this section, we recall some needed ideas regarding
algebraic deformations
(see \cite{gerstenhaber-schack92} for the details), Hochschild
cohomology (see \cite{benson91b}), and crossed products
(see \cite{passman89}). A {\em formal
deformation} of an associative $\C$-algebra $A$ is an algebra
$A[[t]]$ over the formal power series ring $\C[[t]]$ in one
variable, with multiplication defined by
$$u*v = uv +\mu_1(u\otimes v)t + \mu_2(u\otimes v)t^2 + \cdots$$
($u,v\in A$), where the
$\mu_i:A\otimes A\rightarrow A$ are linear maps.
Associativity of $A[[t]]$ imposes constraints on the $\mu_i$. In
particular, the {\em infinitesimal} (or {\em first order}) {\em
deformation} $\mu_1$ must satisfy
\begin{equation}\label{mu1}
\mu_1(u\otimes v)w+\mu_1(uv\otimes w)=\mu_1(u\otimes vw)+u\mu_1(v\otimes w)
\end{equation}
for all $u,v,w\in A$, that is $\mu_1$ is a Hochschild two-cocycle
(a representative of an element in $\HH^2(A)$). Here, as $A$ is an
algebra over a field, its {\em Hochschild cohomology} may be
defined as
$$\HH^*(A):=\Ext^*_{A^e}(A,A),$$
where $A^e:=A\otimes A^{op}$ acts on $A$ by left and
right multiplication. Often this is expressed in terms of the {\em
(acyclic) Hochschild complex}, that is the $A^e$-free resolution
of $A$ given by:
\begin{equation}\label{hc}
\cdots\stackrel{\delta_3}{\longrightarrow}A^{\otimes 4}\stackrel
{\delta_2}{\longrightarrow}A^{\otimes
3}\stackrel{\delta_1}{\longrightarrow}
A^e\stackrel{m}{\longrightarrow}A\longrightarrow 0,
\end{equation}
where $m$ is the multiplication map and
$$\delta_n(a_0\otimes a_1\otimes\cdots\otimes a_{n+1})=
  \sum_{i=0}^n (-1)^i a_0\otimes\cdots\otimes a_ia_{i+1}\otimes\cdots
   \otimes a_{n+1}.$$
Removing the term $A$ from the complex (\ref{hc}) and applying
$\Hom_{A^e}(-,A)$, we obtain the {\em Hochschild (cochain) complex}
\begin{equation}\label{cochain}
0\longrightarrow\Hom_{A^e}(A^e,A)\stackrel{\delta_1^*}{\longrightarrow}
  \Hom_{A^e}(A^{\otimes 3},A)\stackrel{\delta_2^*}{\longrightarrow}
  \Hom_{A^e}(A^{\otimes 4},A)\stackrel{\delta_3^*}{\longrightarrow}\cdots
\end{equation}
Thus $\HH^n(A)=\Ker(\delta_{n+1}^*)/\Ima(\delta_n^*)$. Since
$\Hom_{A^e}(A^{\otimes (n+2)},A)\cong \Hom_{\C}(A^{\otimes n},A)$, we
may identify $\HH^2(A)$ with a subquotient of $\Hom_{\C}(A^{\otimes
2},A)$, the space of infinitesimal deformations of $A$ mentioned
earlier (see (\ref{mu1})).
Obstructions to lifting an infinitesimal deformation $\mu_1$ to a
formal deformation $A[[t]]$ of $A$ lie in $\HH^3(A)$
\cite{gerstenhaber-schack92}.

More generally, if $M$ is any $A$-bimodule (equivalently,
$A^e$-module), then $$\HH^*(A,M):= \Ext^*_{A^e}(A,M).$$

Next we recall the definition of a crossed product ring.
Let $G$ be a finite group
acting by automorphisms on a $\C$-algebra $R$. Let $\alpha:
G\times G\rightarrow R^{\times}$ be a two-cocycle (where
$R^{\times}$ is the group of units of $R$), that is
$$(\rho\cdot\alpha(\sigma,\tau))\alpha(\rho,\sigma\tau)=
\alpha(\rho,\sigma)\alpha(\rho\sigma,\tau)$$
for all $\rho,\sigma,\tau\in G$.
We assume that the image of $\alpha$ is {\em central} in $R$.
(More generally, the image of $\alpha$ need not be central, and
the action of $G$ is twisted by $\alpha$.  
See for example \cite{passman89}.)
In fact we will mainly be interested in two-cocycles $\alpha$ with
image in $\C^{\times}$.
Sometimes we will identify $\alpha$ with its cohomology class in
$\coh^2(G,R^{\times})$. Let $$A:=R\#_{\alpha} G$$ (or
$R\#_{\alpha}\C G$) be the corresponding {\em
crossed product} ring. That is, as a vector space, $A=R\otimes
_{\C}\C G$, where $\C G$ is the group algebra, and the
multiplication is given by
$$(p\otimes\sigma)(q\otimes\tau)=p(\sigma\cdot q)\alpha(\sigma, \tau)
\otimes\sigma\tau$$ for all $p,q\in R$ and $\sigma,\tau\in G$. To
shorten notation, we will write $\overline{\sigma}:= 1\otimes \sigma$,
so $p\overline{\sigma}:=p\otimes \sigma$. (Note that $R$ is subalgebra
of $A$, but $\C G$ is a subalgebra only if $\alpha$ is a coboundary.)
The action of $G$ on $R$ becomes an inner action on $A$, as
$$\sigma\cdot p=\overline{\sigma}p(\overline{\sigma})^{-1}
 \ \ \left(\mbox{where }(\overline{\sigma})^{-1}=\alpha^{-1}
   (\sigma,\sigma^{-1})\overline{\sigma^{-1}}\right)$$
for all $\sigma\in G$, $p\in R$.

We will obtain specific information about Hochschild cohomology
and deformations of $A$ when $R=\C[x,y,z]$ and $G$ is abelian. The
following algebra will be the main example of this paper.

\begin{ex}\label{example}{\em
Let $G=\Z/2\times \Z/2$ be the Klein four-group whose elements
will be denoted $1, a, b, c$. We define an action of $G$ as
automorphisms on $R=\C[x,y,z]$ by
$$\begin{array}{rrrrrrrrr}
  a\cdot x & = & -x , \hspace{.2cm} & a\cdot y & = & y, \hspace{.2cm}
   & a\cdot z & = & -z,\\
  b\cdot x & = & -x, \hspace{.2cm} & b\cdot y & = & -y, \hspace{.2cm}
   & b\cdot z & = & z.
\end{array}
$$
Up to coboundaries, there is exactly one nontrivial two-cocycle
$\alpha: G\times G \ra {\mathbb C}^{\times}$, which we take to be
given by $\alpha(1,\sigma)=1=\alpha
(\sigma,1)=\alpha(\sigma,\sigma)$ for all $\sigma\in G$, and
$$
  \alpha(a,b)={\tt i}=-\alpha(b,a), \ \ \alpha(b,c)=
  {\tt i}=-\alpha(c,b), \ \
  \alpha(c,a)={\tt i}=-\alpha(a,c),
$$
where ${\tt i}=\sqrt{-1}$. Then $A:=R\#_{\alpha} G$ is a crossed
product algebra, as defined above.
}\end{ex}

\section{Hochschild cohomology}
\label{sec:hoch}

We will first state a general result about Hochschild cohomology
of the rings $A=R\#_{\alpha} G$, that is an immediate consequence of a result
of \c Stefan on the Hochschild cohomology of a Hopf Galois extension
$A/R$ \cite{stefan95}. Alternatively, we provide a more
constructive (for our purposes) proof in Section 5
under additional hypotheses. Note that $A$
is an $R$-bimodule under left and right multiplication. The
superscript $G$ in the statement of the following proposition
denotes the subspace of $G$-invariant elements, that is all
elements left unchanged by the action of any $g\in G$.

\begin{prop}\label{stefan}
Let $A=R\#_{\alpha}G$. For each $n\geq 0$,
there is an action of $G$ on $\HH^n(R,A)$ such that
$$
  \HH^n(A)\cong \HH^n(R,A)^G.
$$
\end{prop}

\begin{proof}
By \cite[Cor.\ 3.4]{stefan95} (see also
\cite[Prop.\ 2.3, Prop.\ 2.4 and Thm.\ 3.3]{stefan95}),
there is an action of $G$
on $\HH^n(R,A)$ and a spectral sequence with
$$
  E_2^{m,n} = \coh^m(G, \HH^n(R,A)) \implies
   \HH^{m+n}(A).
$$
As we are working in characteristic 0,
the cohomology of $G$ is concentrated in
degree 0, that is
$$
  \coh^*(G, \HH^n(R,A))=\coh^0(G,\HH^n(R,A))\cong
   \HH^n(R,A)^G.
$$
Therefore $E^{m,n}_2 =E^{m,n}_{\infty}$ and $\HH^n(A)\cong\HH^n(R,A)^G$.
\end{proof}

It is again a result of \c Stefan that a $G$-action on $\HH^n(R,A)$
extending the action inherited from $A$ on
$$\HH^0(R,A)\cong A^R:=\{u\in A\mid up=pu \mbox{ for all } p\in R\}$$
exists and is unique
\cite[Prop.\ 2.4]{stefan95}. Therefore, if we can find such a
group action on $\HH^n(R,A)$, it is necessarily the action to
which Proposition \ref{stefan} refers.

Assuming that the image of $\alpha$ is in $\C^{\times}$, the group
action in the proposition arises from an action of $G$ on the
Hochschild complex (\ref{cochain}) for $R,A$:
\begin{equation}\label{G-cx}
  (\sigma\cdot f)(p_1\otimes \cdots\otimes p_n):=
    \overline{\sigma} f((\sigma^{-1}\cdot p_1)\otimes\cdots \otimes
   (\sigma^{-1}\cdot p_n))(\overline{\sigma})^{-1}
\end{equation}
for all $\sigma\in G$, $f\in\Hom_{R^e}(R^{\otimes n},A)$. However
in case $R$ is a polynomial algebra, we would like to have a group
action on the Koszul complex, as this is the complex with which we
will compute cohomology. We will describe such an action explicitly
in case $R=\C[x,y,z]$ and $G$ is abelian.

Let
$$
  f:=x\otimes 1-1\otimes x, \ \ g:=y\otimes 1-1\otimes y,
   \ \ h:=z\otimes 1-1\otimes z\in R^e.
$$
The {\em Koszul complex} (a free $R^e$-resolution of $R$)
is a complex in which the terms are exterior powers of $R^e$ 
\cite{weibel94}, and it is equivalent to
\begin{equation}\label{koszul}
  0\ra R^e\stackrel{\delta_3}{\longrightarrow}(R^e)^{\oplus 3}\stackrel
  {\delta_2}{\longrightarrow} (R^e)^{\oplus 3}\stackrel{\delta_1}
  {\longrightarrow} R^e \stackrel{m}{\longrightarrow} R\ra 0
\end{equation}
where $m$ is multiplication, and
$$
  \delta_1=(f \ g \ h) , \  \delta_2=\left(\begin{array}{rrr}
                                    -h & 0 & -g\\
                                     0 & -h & f\\
                                     f & g & 0\end{array}\right),
   \mbox{ and } \delta_3=\left(\begin{array}{r} -g\\ f\\
h\end{array}\right).
$$
Dropping the last term $R$ from the complex (\ref{koszul}),
mapping into $A$, and
identifying $\Hom_{R^e} ((R^e)^{\oplus n},A)$ with $A^{\oplus n}$,
we obtain the complex
\begin{equation}\label{complex}
  0\ra A\stackrel{\delta_1^*}{\longrightarrow} A^{\oplus 3}\stackrel
  {\delta_2^*}{\longrightarrow}A^{\oplus 3}\stackrel{\delta_3^*}
  {\longrightarrow} A\ra 0.
\end{equation}
Therefore $\delta_n^*$ is just the transpose of $\delta_n$ in each case,
that is
$$
  \delta_1^*=\left(\begin{array}{r}f\\ g\\ h\end{array}\right) , \
  \delta_2^*=\left(\begin{array}{rrr}
                   -h & 0 & f\\
                    0 & -h & g\\
                    -g & f & 0 \end{array}\right), \mbox{ and }
  \delta_3^*= ( -g \ f \ h).
$$
We have $\HH^n(R,A)=\Ker(\delta_{n+1}^*)/\Ima(\delta_n^*)$.

We will need to compute the action of $G$ on the complex (\ref{complex})
induced by the diagonal action on $R^e$ and its exterior powers.
In the following lemma, we assume that $G$ is abelian, and further
that the action is diagonalized so that $G$ acts by scalar
multiplication on each monomial. We define the symbol $p(\sigma)$ by
\begin{equation}\label{psigma}
\sigma\cdot p=p(\sigma)p \ \ \ (\sigma\in G, p \mbox{ a monomial}),
\end{equation}
that is $p(\sigma)$ denotes the scalar by which $\sigma$ acts on $p$.
The lemma below may be verified by direct computation.

\begin{lem}\label{action}
Let $R=\C[x,y,z]$, let $G$ be an abelian group acting
by scalar multiplication on each of $x,y,z$, and let $\alpha\in\coh^2(G,
\C^{\times})$.
The action of $G$ on $R$ induces the following action
on the cochain complex (\ref{complex}). 
Letting $\sigma\in G$ and
$u,v,w\in A=R\#_{\alpha}G$, the action is given by:

\begin{itemize}

\item[(i)] (degree 0) $\ \ \sigma\cdot u=\overline{\sigma}
u\overline {\sigma}^{-1}$

\item[(ii)] (degree 1)
   $$\sigma\cdot\left(\begin{array}{rrr}u\\ v\\ w\end{array}\right)=
     \left(\begin{array}{rrr}
        x(\sigma^{-1})\overline{\sigma}u\overline{\sigma}^{-1}\\
       y(\sigma^{-1}) \overline{\sigma}v\overline{\sigma}^{-1}\\
        z(\sigma^{-1}) \overline{\sigma}w\overline{\sigma}^{-1}
    \end{array}\right)$$

\item[(iii)] (degree 2)
   $$\sigma\cdot\left(\begin{array}{rrr}u\\ v\\ w\end{array}\right)=
     \left(\begin{array}{rrr}
      xz(\sigma^{-1})\overline{\sigma}u\overline{\sigma}^{-1}\\
      yz(\sigma^{-1}) \overline{\sigma}v\overline{\sigma}^{-1}\\
      xy(\sigma^{-1})\overline{\sigma}w\overline{\sigma}^{-1}
           \end{array}\right)$$

\item[(iv)] (degree 3)
  $\ \ \sigma\cdot u
=xyz(\sigma^{-1})\overline{\sigma}u\overline{\sigma}^{-1}$.
\end{itemize}
\end{lem}

We now have the following algorithm for computing $\HH^*(A)$ under
our hypotheses: Find $\HH^*(R,A)$ by direct calculation using
the cochain complex (\ref{complex}). Then compute the
elements of $\HH^*(R,A)$ invariant
under the $G$-action given in Lemma \ref{action}.
By Proposition \ref{stefan} and the remarks immediately following
its proof, $\HH^*(A)\cong\HH^*(R,A)^G$.

We make one general remark before giving an
example: as an $R$-bimodule, $A$ decomposes into a
direct sum $\oplus_{\sigma\in G} R\overline{\sigma}$, and so
$$\HH^*(R,A)\cong \bigoplus_{\sigma\in G}\HH^*(R,R\overline{\sigma}).$$
If $G$ is abelian, the action of $G$ on cohomology preserves this
decomposition, and by Proposition \ref{stefan} we have an
additive decomposition
$$\HH^*(A)\cong \bigoplus_{\sigma\in G}\HH^*(R,R\overline{\sigma})^G.$$
This decomposition of Hochschild cohomology is illustrated in the
following example.

\begin{ex}\label{hh*}{\em
Let $A=R\#_{\alpha}(\Z/2\times \Z/2)$ as in Example \ref{example},
where $\alpha$ may be either trivial or nontrivial. Note that the
$R$-bimodule structure of $A$ is the same in either of the cases
that $\alpha$ is trivial or nontrivial, so $\HH^*(R,A)$ is the
same in either case. We find $\HH^0(R,A)\cong A^R=R\cong\HH^0(R)$,
$\HH^n(R,A)=0$ ($n>3$), and
\begin{itemize}
\item[(i)] $\HH^1(R,A)\cong R^{\oplus 3}\cong \HH^1(R)$
\item[(ii)]
   $\displaystyle{\HH^2(R,A)=\left\{\left(\begin{array}{c}
     p_1+q_1\overline{a}\\ p_2+q_2\overline{c}\\ p_3+q_3
    \overline{b}\end{array}\right)\mid
     p_i\in R, q_1\in  \C [y], q_2\in \C [x], q_3\in \C [z]\right\}}$

  $\hspace{1cm}=R^{\oplus 3}\oplus \C[y]\overline{a}\oplus
\C[z]\overline{b}\oplus
  \C[x]\overline{c}$
\item[(iii)] $\HH^3(R,A)=\{p_1+p_2\overline{a}+p_3\overline{b}
+p_4
  \overline{c}\mid p_1\in R, p_2\in \C[y],p_3\in \C[z],p_4\in \C[x]\}$

  $\hspace{1cm}=
  R\oplus \C[y]\overline{a}\oplus \C[z]\overline{b}\oplus \C[x]\overline{c}$.
\end{itemize}
Applying Lemma \ref{action}, we have $\HH^0(A)\cong
\C[x^2,y^2,z^2,xyz]=Z(A)$, the center of $A$ (as expected), $\HH^n(A)=0$
($n>3$), and
\begin{itemize}
\item[(iv)] In either the case that $\alpha$ is trivial or that
$\alpha$
  is nontrivial,

\noindent
  $\displaystyle{\HH^1(A)=
\left\{\left(\begin{array}{c}p_1\\p_2\\p_3\end{array}\right)
     \mid \begin{array}{l}p_1\in x\C[x^2,y^2,z^2]+yz\C[x^2,y^2,z^2]\\
         p_2\in y\C[x^2,y^2,z^2]+xz\C[x^2,y^2,z^2]\\
         p_3\in z\C[x^2,y^2,z^2]+xy\C[x^2,y^2,z^2]
     \end{array}\right\}}$
\item[(v)] $\HH^2(A)=\left\{\left(\begin{array}{c}p_1+q_1\overline{a}\\
  p_2+q_2\overline{c}\\p_3+q_3\overline{b}\end{array}\right)\right\}$,
where $p_1\in y\C [x^2,y^2,z^2] + xz\C [x^2,y^2,z^2]$,\\ $p_2\in x\C
[x^2,y^2,z^2]+yz \C [x^2,y^2,z^2]$, and $p_3\in z\C
[x^2,y^2,z^2]+xy\C [x^2,y^2,z^2]$. If $\alpha$ is trivial, we take
$q_1\in y\C [y^2]$, $q_2\in x\C [x^2]$, and $q_3\in z\C [z^2]$. If
$\alpha$ is nontrivial, we take $q_1\in \C [y^2]$, $q_2\in \C [x^2]$,
and $q_3\in \C [z^2]$.

\item[(vi)]
$\HH^3(A)=\{p_1+p_2\overline{a}+p_3\overline{b}+p_4\overline{c}\}$,
where $p_1\in Z(A)$. If $\alpha$ is trivial, we take $p_2\in
\C[y^2]$, $p_3\in \C[z^2]$, and $p_4\in\C[x^2]$. If $\alpha$ is
nontrivial, we take $p_2\in y\C[y^2]$, $p_3\in z\C[z^2]$, and $p_4\in
x\C[x^2]\}$.
\end{itemize}
}\end{ex}

\section{Infinitesimal deformations}
\label{sec:infidef}
Under some hypotheses on $A=R\#_{\alpha}G$,
we find in this section an explicit formula for the infinitesimal deformation
$\mu_1:A\otimes A\rightarrow A$ corresponding to a given
element of $\HH^2(A)$.
We assume that $\HH^2(A)$ has
been computed via Proposition \ref{stefan}, using an
$R^e$-projective resolution of $R$ that itself carries an action of $G$.
(See for example Lemma \ref{action}.)
Since our elements of
$\HH^2(A)$ are given as $G$-invariant elements of $\HH^2(R,A)$, we
will therefore need to relate the Hochschild complex for $A$ to this
resolution of $R$. With this in mind, we will now give a more constructive
(for our purposes)
proof of Proposition \ref{stefan}, using the action (\ref{G-cx}) of
$G$ on the Hochschild complex.

We assume now that the image of $\alpha$ is contained in $\C^{\times}$.
Let
$$ \Delta := \bigoplus_{\sigma\in G} R\overline{\sigma} \otimes R
\overline{\sigma}^{-1},$$ a subalgebra of $A^e=A\otimes A^{op}$
containing $R^e$. Then $R$ is a $\Delta$-module and
$A\cong (A^e)\otimes_{\Delta} R$ as an $A^e$-module (see
\cite[Lemma 3.3]{boisen92}).
Note that $A^e=\oplus_{\sigma\in
G}\Delta(\overline{\sigma} \otimes 1)$ is a free $\Delta$-module.
Therefore by the Eckmann-Shapiro Lemma \cite[Corollary
2.8.4]{benson91a},
$$ \Ext^n_{A^e}(A,A)\cong \Ext^n_{A^e}\left((A^e)
\otimes_{\Delta}R, A\right)\cong \Ext^n_{\Delta}(R,A).$$ Once we
show that $\Ext^n_{\Delta}(R,A)\cong
\left(\Ext^n_{R^e}(R,A)\right) ^G$, we will have proved
Proposition \ref{stefan}.

Let $ \cdots \stackrel{\delta_3}{\longrightarrow} P_2
\stackrel{\delta_2}{\longrightarrow}
P_1\stackrel{\delta_1}{\longrightarrow}
P_0\stackrel{\delta_0}{\longrightarrow} R \longrightarrow 0 $ be a
$\Delta$-projective resolution of $R$. By restriction, it is an
$R^e$-projective resolution of $R$, and there is an inclusion
$\Hom_{\Delta}(P_n,A)\subset\Hom_{R^e}(P_n,A)$ for each $n$.
Moreover, there is an action of $G$ on the complex given by
$\sigma\cdot p=\overline{\sigma}p\overline{\sigma}^{-1}$ (as the
image of $\alpha$ is in $\C^{\times}$), and a
corresponding action on the cochain complex given by $(\sigma\cdot
f)(p)=\overline{\sigma}f((\overline{\sigma})^{-1}p
\overline{\sigma})(\overline{\sigma})^{-1}$ for all $\sigma\in G,
p\in P_n,$ and $f\in\Hom_{R^e}(P_n,A)$. Therefore
$\Hom_{\Delta}(P_n,A)=\Hom_{R^e}(P_n,A)^G$. As $|G|$ is invertible
in $\C$, the subspace of $G$-invariant elements is just the image
of the trace map $\frac{1}{|G|}\sum_{\sigma\in G}\sigma\cdot$, and
so $(\Ker(\delta_{n+1}^*)/\Ima(\delta_n^*))^G\cong
\Ker(\delta_{n+1}^*)^G/ \Ima(\delta_n^*)^G$. Further, if
$f=\delta_n^*(f')$ is $G$-invariant, then
$f=\frac{1}{|G|}\sum_{\sigma\in G}\sigma\cdot\delta_n^*(f') =
\delta_n^*(\frac{1}{|G|}\sum_{\sigma\in G}\sigma\cdot f')$ is in
the image of $\delta_n^*$ restricted to $G$-invariant
homomorphisms, that is $\Delta$-homomorphisms. Therefore
$\Ext_{\Delta}^n(R,A)\cong \left(\Ext^n_{R^e}(R,A)\right)^G$,
which proves Proposition \ref{stefan} under the assumption that
$\alpha\in \coh^2(G,\C^{\times})$.

We will find useful the following $\Delta$-projective resolution of
$R$.
For each $n\geq 0$, let
$$\Delta_n:=\left\{\sum p_0\overline{\sigma_0}\otimes \cdots \otimes
p_{n+1}\overline{\sigma_{n+1}}\mid p_i\in R, \sigma_i\in G \mbox{ and }
\sigma_0\cdots \sigma_{n+1}=1\right\},$$
a $\Delta$-submodule of $A^{\otimes (n+2)}$.
Thus $\Delta_0=\Delta$, and each $\Delta_n$ is a projective
$\Delta$-module. Consider the complex
\begin{equation}\label{Delta}
  \cdots \stackrel{\delta_3}{\longrightarrow}
  \Delta_2 \stackrel{\delta_2}{\longrightarrow}
  \Delta_1 \stackrel{\delta_1}{\longrightarrow}\Delta_0\stackrel
  {m}{\longrightarrow} R \longrightarrow 0
\end{equation}
where $m$ denotes multiplication and the maps $\delta_n$ are
restrictions of the standard maps from the Hochschild complex
(\ref{hc}). The above complex (\ref{Delta}) is exact, as there is
a chain contraction $s_n:\Delta_{n-1}\rightarrow \Delta_n$ given
by
$$s_n(u_0\otimes \cdots\otimes u_n)=u_0\otimes \cdots \otimes u_n\otimes
1.$$
Therefore (\ref{Delta}) is a $\Delta$-projective resolution of $R$.

The Hochschild complex (\ref{hc}) is just the complex
(\ref{Delta}) induced from $\Delta$ to $A^e$. Thus the isomorphism
$\Ext^n_{A^e}(A,A)\stackrel{\sim} {\longrightarrow}
\Ext^n_{\Delta}(R,A)$ is given at the chain level simply by
restricting maps from $\Hom_{A^e}(A^{\otimes (n+2)},A)$ to
$\Hom_{\Delta}(\Delta_n,A)$. As discussed above, the isomorphism
$\Ext^n_{\Delta}
(R,A)\stackrel{\sim}{\longrightarrow}\Ext^n_{R^e}(R,A)^G$ is given
at the chain level by inclusion of $\Hom_{\Delta}(\Delta_n, A)$
into $\Hom_{R^e}(\Delta_n,A)$. Now there is a map from
(\ref{Delta}) to the Hochschild complex for $R$, as they are both
$R^e$-projective resolutions of $R$. Under the assumption that
$\alpha\in \coh^2(G,\C^{\times})$, this is given by
\begin{equation*}
  p_0\overline{\sigma_0}\otimes p_1\overline{\sigma_1}\otimes p_2
  \overline{\sigma_2}\otimes\cdots\otimes p_{n+1}\overline{\sigma_{n+1}}
  \mapsto\hspace{3in}
\end{equation*}
\begin{equation}\label{mapto}
\hspace{.2in}\left(\prod_{i=0}^n \alpha(\sigma_0\cdots \sigma_i,
\sigma_{i+1})\right)p_0\otimes \sigma_0\cdot p_1\otimes
(\sigma_0\sigma_1)\cdot p_2\otimes
\cdots\otimes(\sigma_0\sigma_1\cdots \sigma_n)\cdot p_{n+1}.
\end{equation}
In particular, this map is the identity on the submodule
$R^{\otimes (n+2)}$ of $\Delta_n$.

Now suppose that
\begin{equation}\label{P}
\cdots P_2\rightarrow P_1\rightarrow P_0\rightarrow R\rightarrow 0
\end{equation}
is any other $R^e$-projective resolution of $R$ carrying an action of $G$.
Let $$\psi_n:R^{\otimes (n+2)}\rightarrow P_n \
(n\geq 0)$$ be $R^e$-homomorphisms giving a map of chain complexes
from the Hochschild complex (\ref{hc}) for $R$ to the above complex
(\ref{P}). The next
result is an explicit formula for the infinitesimal deformation
$\mu_1:A\otimes A\rightarrow A$ corresponding to an element of
$\HH^2(A)\cong \HH^2(R,A)^G$ expressed as a cochain in terms of
(\ref{P}).

\begin{thm} \label{mu}
Let $A=R\#_{\alpha}G$ with $\alpha\in \coh^2(G,\C^{\times})$.
Let $f: P_n\rightarrow A$ be a function representing an element of
$\HH^n(R,A)^G\cong \HH^n(A)$ expressed in terms of the complex (\ref{P}).
The corresponding function
$\tilde{f}\in\Hom_{\C}(A^{\otimes n},A)\cong\Hom_{A^e}(A^{\otimes (n+2)},
A)$ of the Hochschild complex (\ref{cochain}) is given by
$$\tilde{f}(p_1\overline{\sigma_1}\otimes\cdots\otimes
p_n\overline{\sigma_n})=
((f\circ \psi_n)(1\otimes p_1\otimes\sigma_1 \cdot
p_2\otimes\cdots\otimes (\sigma_1\cdots\sigma_{n-1})\cdot
p_n\otimes 1))\overline{\sigma_1}\cdots \overline{\sigma_n}.$$ In
particular, if $n=2$, we obtain the infinitesimal deformation
$\mu_1:A\otimes A\rightarrow A$,
$$\mu_1 (p\overline{\sigma}\otimes q\overline{\tau})=
  ((f\circ \psi_2)(1\otimes p\otimes\sigma\cdot q \otimes 1))
  \overline{\sigma}\cdot\overline{\tau}.$$
\end{thm}

\begin{proof}
Identifying $\Hom_{\C}(A^{\otimes n},A)$ with $\Hom_{A^e}(A^{\otimes
(n+2)},A)$ and applying the map (\ref{mapto}), we have
\begin{eqnarray*}
  && \hspace{-1cm}
  \tilde{f}(1\otimes p_1\overline{\sigma_1}\otimes\cdots\otimes p_n
  \overline{\sigma_n}\otimes 1)\\
&=&\alpha^{-1}(\sigma_1\cdots\sigma_n,\sigma_n^{-1}\cdots\sigma_1^{-1})
   \tilde{f}\left(1\otimes p_1\overline{\sigma_1}\otimes \cdots\otimes p_n
\overline{\sigma_n}\otimes\overline{\sigma_n^{-1}\cdots\sigma_1^{-1}}\right)
  \overline{\sigma_1\cdots\sigma_n}\\
   &=& \prod_{i=1}^{n-1}\alpha(\sigma_1\cdots\sigma_i,\sigma_{i+1})((f\circ
   \psi_n)(1\otimes p_1\otimes \sigma_1\cdot p_2\otimes\cdots\otimes
    (\sigma_1\cdots\sigma_{n-1})\cdot p_n)\otimes 1))\overline{\sigma_1
  \cdots\sigma_n}\\
   &=& ((f\circ \psi_n)(1\otimes p_1\otimes\sigma_1\cdot p_2\otimes\cdots
   \otimes (\sigma_1\cdots\sigma_{n-1})\cdot p_n\otimes
1)\overline{\sigma_1}
   \cdots \overline{\sigma_n}.
\end{eqnarray*}\end{proof}

We will use Theorem \ref{mu} to calculate the infinitesimal
deformations for Example \ref{example}, that is for $A=\C[x,y,z]\#_{\alpha}
(\Z/2\times \Z/2)$ (whose Hochschild cohomology was given in
Example \ref{hh*}).
Let $R=\C[x,y,z]$ for the
rest of this section. There are $R^e$-homomorphisms giving a map
from the Hochschild complex (\ref{hc}) to the Koszul complex
(\ref{koszul}) for $R$,
$$\begin{array}{ccccccccccc}
\cdots\stackrel{\delta_4}{\longrightarrow} & R^{\otimes 5} &
\stackrel {\delta_3}{\longrightarrow} & R^{\otimes 4}&
\stackrel{\delta_2} {\longrightarrow} & R^{\otimes 3}&
\stackrel{\delta_1}{\longrightarrow} &
R^{\otimes 2} & \stackrel{m}{\longrightarrow} & R & \longrightarrow 0\\
&\downarrow\psi_3 && \hspace{.2in}\downarrow\psi_2 &&
\hspace{.2in}\downarrow\psi_1
&& || && || &\\
0\longrightarrow & R^e& \stackrel{\delta_3}{\longrightarrow} &
(R^e)^{\oplus 3} & \stackrel{\delta_2}{\longrightarrow} &
(R^e)^{\oplus 3} & \stackrel{\delta_1}{\longrightarrow} & R^e &
\stackrel{m}{\longrightarrow} & R & \longrightarrow 0.
\end{array}$$
A computation shows that such maps may be defined as follows
(cf.\ \cite{halbout01}):
$$\psi_1(1\otimes x^iy^jz^k\otimes 1) = \left(\begin{array}{c}
\displaystyle{\sum_{\ell =1}^i} x^{i-\ell}y^jz^k\otimes x^{\ell -1}\\
\displaystyle{\sum_{\ell = 1}^j} y^{j-\ell}z^k\otimes y^{\ell -1}x^i\\
\displaystyle{\sum_{\ell =1}^k} z^{k-\ell}\otimes x^iy^jz^{\ell
-1}\end{array}\right)$$
$$\psi_2(1\otimes x^iy^jz^k\otimes x^ry^sz^t\otimes 1)=\hspace{3.5in}$$
$$\left(\begin{array}{c}\displaystyle{\sum_{m=1}^t\sum_{\ell =1}^i}
x^{i-\ell}
y^{j+s}z^{k+t-m}\otimes x^{r+\ell -1}z^{m-1}\\
\displaystyle{\sum_{m=1}^t\sum_{\ell
=1}^{j+s}}y^{j+s-\ell}z^{k+t-m} \otimes x^{i+r}y^{\ell -1}z^{m-1}
- \displaystyle{\sum_{m=1}^t\sum_{\ell =1}^s} x^iy^{j+s-\ell}
z^{k+t-m}\otimes x^ry^{\ell -1}z^{m-1}\\
\displaystyle{\sum_{m=1}^s\sum_{\ell =1}^i}
x^{i-\ell}y^{j+s-m}z^k\otimes x^{r+\ell -1} y^{m-1}z^t
\end{array}\right)$$
$$\psi_3(1\otimes x^iy^jz^k \otimes x^ry^sz^t\otimes x^uy^vz^w\otimes 1)=
\hspace{3.5in}$$
$$-\sum_{n=1}^v\sum_{m=1}^t\sum_{\ell =1}^i x^{i-\ell}y^{j+s+n-1}
z^{k+t-m}\otimes x^{r+u+\ell -1}y^{v-n}z^{w+m-1}.$$

\begin{ex}\label{infinitesimal}{\em
Let $A=\C[x,y,z]\#_{\alpha}(\Z/2\times\Z/2)$ as in Example
\ref{example}, where $\alpha$ may be either trivial or nontrivial.
For any integer $i$, let $\overline{i}$ be its reduction modulo 2,
that is $\overline{i}=0$ if $i$ is even and $\overline{i}=1$ if
$i$ is odd. Direct computation yields the following:
the infinitesimal deformation $\mu_1: A\otimes
A\rightarrow A$ given by Theorem \ref{mu} and map $\psi_2$ above,
and corresponding to the element $(p_1,p_2,p_3; q_1\overline{a};
q_2\overline{c}; q_3\overline{b})$ of Example \ref{hh*}(v) is
$$\begin{array}{ll}
\mu_1(x^iy^jz^k \overline{\sigma} \otimes x^{\ell}y^mz^n \overline{\tau})  &
\\
\hspace{2cm} = x^{\ell}y^mz^n(\sigma) &\hspace{-.3cm}\left\{ (x^{i+\ell
-1}y^{j+m}z^{k+n-1})\left(inp_1 + \overline{in} (-1)^{\ell} q_1
\overline{a}\right)\right.\\
& \hspace{.4cm}
+(x^{i+\ell}y^{j+m-1}z^{k+n-1})\left(jn p_2 + \overline{jn}(-1)^m q_2
\overline{c}\right)\\
 & \left.\hspace{.4cm}+(x^{i+\ell -1}y^{j+m-1}z^{k+n})\left(imp_3 +
\overline{im}(-1)^{\ell} q_3 \overline{b}\right)\right\}\overline{\sigma}\cdot
\overline{\tau}.
\end{array}$$
Here the scalar $x^{\ell}y^mz^n(\sigma)$ is defined by (\ref{psigma}), and
a negative power of $x,y$ or $z$ indicates that the term is 0.
}\end{ex}

\section{Formal deformations}
\label{sec:formal}

We will give a universal deformation formula (definition below)
that will apply in particular to $A=\C[x,y,z]\#_{\alpha}(\Z/2\times \Z/2)$,
producing a formal
deformation that lifts some of its infinitesimal deformations
in case $\alpha$ is nontrivial. Our formula was inspired by an
unpublished result of the second author and Zhang, and uses results in
their paper \cite{giaquinto-zhang98}.

A {\em universal deformation
formula} based on a bialgebra $B$ is an element $F\in (B\otimes
B)[[t]]$ satisfying the equations
\begin{multline}\label{udf-def}
\hspace{3cm}(\epsilon\otimes\id)(F) = 1\otimes 1 = (\id\otimes\epsilon)(F)\\
\mbox{and } \ [(\Delta\otimes\id)(F)]\cdot (F\otimes 1)
= [(\id\otimes\Delta)(F)]\cdot (1\otimes F),\hspace{2cm}
\end{multline}
where $\id$ is the identity map. The virtue of such
a formula is that if $S$ is any $B$-module algebra, then $F$
provides a formal deformation of $S$ \cite{giaquinto-zhang98}.
Specifically, if $x,y\in S$,
then the deformed product is given by $x*y=(m\circ F)(x\otimes
y)$, where $m: S\otimes S\rightarrow S$ is the ordinary
multiplication.

Let
$H_1$ be the associative $\C$-algebra generated by the elements
$D_1$, $D_1'$, and $\beta_1$, subject to the relations
$$D_1^2=0=(D_1')^2 , \ \ D_1D_1'=D_1'D_1, \ \ D_1\beta_1=-\beta_1D_1,
\ \ D_1'\beta_1=-\beta_1D_1',  \mbox{ and } \beta_1^2=1.$$ Then
$H_1$ is a bialgebra with comultiplication determined by
$$\Delta(D_1)=D_1\otimes \beta_1 + 1\otimes D_1, \ \
\Delta(D_1')=D_1'\otimes 1 + \beta_1\otimes D_1', \ \
\Delta(\beta_1)=\beta_1\otimes \beta_1,$$ and counit
$\epsilon(D_1)=0=\epsilon(D_1')$, $\epsilon(\beta_1)=1$. In fact
$H_1$ is a quotient of the Drinfel'd double of a Taft algebra (see
\cite[Lemma 4.4]{montgomery-schneider01}). Let
$H_1[[t]]=H_1\otimes_{\C}\C[[t]]$ be the Hopf algebra $H_1$ with
coefficients extended to the formal power series ring $\C[[t]]$.
Take $H_2$ and $H_3$ to be two more copies of
the Hopf algebra $H_1$ (with appropriately changed indices), and
$H=H_1\otimes H_2\otimes H_3$, the tensor product Hopf algebra.

\begin{lem}\label{udf}
The element $F_i=1\otimes 1 + t D_i\otimes D_i'$ of
$H_i[[t]\otimes H_i[[t]]$ ($i=1,2,3$) is a universal deformation formula.
Consequently, the product $F=F_1F_2F_3$ is a universal deformation
formula, based on $H[[t]]\otimes H[[t]]$.
\end{lem}

\begin{proof}
It is straightforward to check the equations (\ref{udf-def}) for each $F_i$.
The second statement of the lemma follows, as
the product of universal deformation formulas is again a universal
deformation formula, based on the tensor product of the
bialgebras.
\end{proof}

\begin{ex}\label{example5}{\em
Let $A=\C[x,y,z]\#_{\alpha}(\Z/2\times\Z/2)$, as in Example
\ref{example}. We assume the cocycle $\alpha\in \coh^2(\Z/2\times
\Z/2,\C^{\times})$ is {\em nontrivial}. We claim that $A$ is an
$H$-module algebra, that is the $D_i$ and $D_i'$ act as skew
derivations, and the $\beta_i$ act as automorphisms. Let
$q_1\in\C[y^2]$, $q_2\in\C[x^2]$, and $q_3\in\C[z^2]$ (see Example
\ref{hh*}(v)). Define, for all $i,j,k\in \Z^{\geq 0}$ and
$\sigma\in \Z/2\times \Z/2$,
$$D_1(x^iy^jz^k\overline{\sigma}) =  \overline{i}
x(\sigma)x^{i-1}y^jz^k\overline{\sigma}, \ \ \
D_1'(x^iy^jz^k\overline{\sigma}) =
(-1)^i\overline{k}x^iy^jz^{k-1}q_1 \overline{a}\cdot
\overline{\sigma},$$
$$D_2(x^iy^jz^k\overline{\sigma}) =  \overline{j}
z(\sigma) x^iy^{j-1}z^kq_2\overline{c} \cdot \overline{\sigma}, \
\ \ D_2'(x^iy^jz^k\overline{\sigma}) = \overline{k}x^iy^jz^{k-1}
\overline{\sigma},$$
$$D_3(x^iy^jz^k\overline{\sigma}) =  \overline{i}
y(\sigma) x^{i-1} y^jz^kq_3\overline{b}\cdot \overline{\sigma} , \
\ \ D_3'(x^iy^jz^k\overline{\sigma}) = \overline{j}
x^iy^{j-1}z^k\overline{\sigma},$$
$$\beta_1(x^iy^jz^k\overline{\sigma}) =  (-1)^i
x(\sigma)x^iy^jz^k\overline{\sigma}, $$
$$\beta_2(x^iy^jz^k\overline{\sigma}) = (-1)^k
z(\sigma)x^iy^jz^k\overline{\sigma},$$
$$\beta_3(x^iy^jz^k\overline{\sigma})=(-1)^j
y(\sigma)x^iy^jz^k\overline{\sigma},$$
where the scalars $x(\sigma)$, $y(\sigma)$, and $z(\sigma)$ are
defined in (\ref{psigma}).
Then $A$ is a left
$H$-module algebra under the above operations: from their
definitions, it is clear that $\beta_1$, $\beta_2$, and $\beta_3$
are automorphisms on $A$, and that they commute with each other.
It is straightforward to verify that the $D_i$, $D_i'$ are skew
derivations on $A$ with respect to the automorphisms $\beta_i$, in
accordance with their coproducts in $H\otimes H$. From their
definitions we have $D_i^2=0=(D_i')^2$ and $\beta_i^2=\id$. Again
it may be verified directly that $D_i\beta_i=-\beta_iD_i$, and
$D_i'\beta_i=-\beta_iD_i'$ as operators on $A$ ($i=1,2,3$), and
that all other pairs of these operators commute.

Now set $p_1=p_2=p_3=0$ in Example \ref{infinitesimal}. We claim
that the corresponding infinitesimal deformation $\mu_1$ of $A$
lifts to a formal deformation of $A$ over the power series ring
$\C[[t]]$: by \cite[Thm.\ 1.3 and Defn.\ 1.11]{giaquinto-zhang98}
and Lemma \ref{udf}, the universal deformation formula $F$ gives a
formal deformation of $A$.  From their definitions, as operators
on $A$, we see that
$$D_1D_3 \mid _{A}=0, \ \ D_2D_3' \mid _{A}=0, \mbox{ and }
D_1'D_2' \mid _{A}=0.$$ Therefore we may write the universal
deformation formula $F$, as an operator on $A\otimes A$, as
\begin{equation}\label{F}
F | _{A\otimes A}=1\otimes 1 + t (D_1\otimes D_1' + D_2\otimes
 D_2' + D_3\otimes D_3') + t ^2 D_2D_3\otimes D_2'D_3'.
\end{equation}
It may be verified that $D_1\otimes D_1' + D_2\otimes D_2' +
D_3\otimes D_3'$, composed with multiplication in $A$, yields the
infinitesimal deformation $\mu_1$ given in Example
\ref{infinitesimal}.

In case the functions $p_i$ are {\em not} all 0, we recover the
infinitesimal deformation $\mu_1$ of Example \ref{infinitesimal}
from the following skew derivations $d_i,d_i'$ and automorphisms
$\gamma_i$:
$$d_1(x^iy^jz^k\sigma)=iz(\sigma){z}x^{i-1}y^jz^kp_1\sigma, \ \ \
d_1'(x^iy^jz^k\sigma)=kx^iy^jz^{k-1}\sigma,$$
$$d_2(x^iy^jz^k)=jz(\sigma){z}x^iy^{j-1}z^kp_2\sigma, \ \ \
d_2'=d_1',$$
$$d_3(x^iy^jz^k\sigma)=iy(\sigma)x^{i-1}y^jz^kp_3\sigma, \ \ \
d_3'(x^iy^jz^k\sigma)=jx^iy^{j-1}z^k\sigma ,$$
$$\gamma_1(x^iy^jz^k\sigma)=z(\sigma)x^iy^jz^k\sigma, \ \ \
\gamma_2=\gamma_1, \ \ \ \gamma_3(x^iy^jz^k\sigma)=y(\sigma)
x^iy^jz^k\sigma.$$ The functions $d_1\otimes d_1'$, $d_2\otimes
d_2'$, and $d_3\otimes d_3'$ on $A\otimes A$ produce the
infinitesimal deformations $\mu_1$ of Example \ref{infinitesimal}
corresponding to $p_1$, $p_2$, and $p_3$, respectively. However,
it may be checked for example that $d_1d_1'\neq d_1'd_1$ in
general, so $A$ is not an $H_1$-module algebra under the actions
of $d_1,d_1'$ and $\gamma_1$. Therefore Lemma \ref{udf} does not
produce a formal deformation of $A$ lifting the infinitesimal
$\mu_1$ in this case. }\end{ex}

\begin{rem}{\em
In case $A=\C[x,y,z]\#_{\alpha}(\Z/2\times\Z/2)$ with $\alpha$
{\em trivial}, $A$ is an $H_i$-module algebra for each $i$, under
the operators $D_i,D_i'$ and $\beta_i$ defined above, but with
$q_1\in y\C[y^2], q_2\in x\C[x^2]$ and $q_3\in z\C[z^2]$ (see
Example \ref{hh*}(v)). So for example, if $q_1\in y\C[y^2]$ and
$q_2=q_3=p_1=p_2=p_3=0$, the universal deformation formula
$F=1\otimes 1 + tD_1\otimes D_1'$ of Lemma \ref{udf} yields a
formal deformation of $A$ lifting the infinitesimal deformation
$\mu_1$ of Example \ref{infinitesimal}. However, $D_1D_2\neq
D_2D_1$ in general in this case, and so $A$ is not an $H$-module
algebra, and the universal deformation formula of Lemma \ref{udf}
based on $H[[t]]$ does not apply to $A$ in case $q_1q_2\neq 0$.
}\end{rem}

\section{The center of the deformed algebra}

In this section, we consider only the example
$A=\C[x,y,z]\#_{\alpha} (\Z/2\times\Z/2)$, where $\alpha$ is {\em
nontrivial}. Let $A_F$ denote the deformed algebra over $\C[[t]]$
given in Example \ref{example5}, in the case $p_1=p_2=p_3=0$. The
center of the original algebra $A$ is generated by $x^2$, $y^2$,
$z^2$, and $xyz$. Inspection of  the universal deformation formula
$F$ given in (\ref{F}) shows that the images of $x^2$, $y^2$, and
$z^2$ remain central in $A_F$, but that this is not necessarily
the case for $xyz$. In the next lemma, we adjust the element $xyz$
so that it is central in $A_F$. The resulting
elements generate the new center.

\begin{lem} Let $F$ be the universal deformation formula given in
(\ref{F}) above. Let $Z$ be the ${\mathbb C}[[t]]$-subalgebra of
$A_F$ generated by $x^2$, $y^2$, $z^2$, and
$$w=xyz+\frac{1}{2}t(yq_1\overline{a}+xq_2
\overline{c}+zq_3\overline{b}).$$ The center of $A_F$ is the
completion $\widehat{Z}$ of $Z$ with respect to the ideal $(t)$.
\end{lem}

\begin{proof}
The element $w$ will be central in $A_F$ once we see
that $w*x^iy^jz^k\overline{\sigma}=x^iy^jz^k
\overline{\sigma} * w$ for any $\sigma\in G$ and nonnegative
integers $i,j,k$. In either product, the resulting term in $A$ is
$x^{i+1}y^{j+1}z^{k+1}\overline{\sigma}$. Lengthy calculations show that
the coefficients of $t$, $t^2$, and $t^3$ in the products $w *
x^iy^jz^k\overline{\sigma}$ and $x^iy^jz^k\overline{\sigma} *w$ are equal.
Clearly this implies that $\widehat{Z}$ is contained in the center of $A_F$.

Now suppose that
\begin{equation}\label{v}
v=\sum_{i=0}^{\infty}t^iv_i
\end{equation}
is in the center of $A_F$.  Then $v_0$ must be in the center of
$A$, and so is equal to a polynomial in $x^2$, $y^2$, $z^2$, and
$xyz$, say $v_0=P_0(x^2,y^2,z^2,xyz)$.  Let
$v_0'=P_0(x^2,y^2,z^2,w)$ as an element of $A_F$, known to be
central by our previous argument. Note that $v_0'=v_0\mod (t)$.
Write
$$v=v_0'+\left((v_0-v_0')+\sum_{i=1}^{\infty}t^iv_i\right),$$
a sum of an element in $Z$ and an element with constant term 0.
Again, the coefficient of $t$ in the expression
$(v_0-v_0')+\sum_{i=1}^{\infty} t^iv_i$ must  be central in $A$,
and so must be a polynomial in $x^2$, $y^2$, $z^2$, and $xyz$, say
$P_1(x^2,y^2,z^2,xyz)$. Let $v_1'=P_1(x^2,y^2,z^2,w)$, and note
that $v_1'$ is equivalent to the coefficient of $t$ in
$(v_0-v_0')+\sum_{i=1}^{\infty}t^iv_i$, modulo $(t^2)$. By
induction, for any positive integer $n$, we may find
$v_0',v_1',\ldots, v_{n-1}'\in \Z$ with
$$v=v_0'+tv_1'+\cdots + t^{n-1}v_{n-1}' \mod (t^n).$$
As $\widehat{Z}:=\displaystyle{\lim_{\longleftarrow}} \ Z/(t^n)$, we may
identify $v$ with an element of $\widehat{Z}$.
\end{proof}

Finally, for homogeneous $q_1,q_2$ and $q_3$, we will find a
presentation of the central subalgebra $Z$ of $A_F$ generated by
$x^2$, $y^2$, $z^2$, and $w$.  To this end, we calculate
$$w*w=x^2y^2z^2+\frac{1}{4}t^2(y^2q_1^2+x^2q_2^2+z^2q_3^2)
-\frac{1}{4}{\tt i}t^3 q_1q_2q_3,$$
where ${\tt i}=\sqrt{-1}$.
Letting $\hat{x}=x^2$, $\hat{y}=y^2$,
$\hat{z}=z^2$, and $\hat{w}=w$ in $A_F$, and choosing $q_1=\hat{y}^j/2$,
$q_2=\hat{x}^i/2$, and $q_3=\hat{z}^k/2$ for nonnegative integers $i,j,k$, we
have a single relation in $Z$,
\begin{equation}\label{relation}
\hat{w}^2=\hat{x}\hat{y}\hat{z} +t^2(\hat{x}^{2i+1}+\hat{y}^{2j+1}+
\hat{z}^{2k+1})-2{\tt i}t^3 \hat{x}^i\hat{y}^j\hat{z}^k.
\end{equation}
Thus $Z$ is generated by $\hat{x},\hat{y},\hat{z}$, and $\hat{w}$,
subject to the relation (\ref{relation}).

A few comments are in order here:
\begin{enumerate}
\item The deformed defining equation~(\ref{relation}) of the center
$Z$ has no terms in $t^1$.  Therefore, to first order $Z$ does not
deform.  This explains our statement in the introduction that the map
$\Def_{R\#_\alpha G} \ra \Def_{R^G}$ is totally ramified at 0.
\item Unlike the compact situation of Vafa and Witten~\cite{VafWit},
we never obtain ordinary double points as our singularities.  In this
respect the local picture is more similar to the second example
studied in [loc.\ cit.].
\item The deformation smooths out the initial three curves of
singularities, leaving only a singularity at the origin, as can be
verified by direct calculation.
\end{enumerate}

\end{document}